\documentclass[11pt]{article} 

\usepackage{times}
\usepackage{amssymb}
\usepackage{prooftree}
\usepackage{amsmath}
\usepackage{latexsym}
\usepackage{stmaryrd}

\newtheorem{definition}{Definition}[section]

\newtheorem{remark}[definition]{Remark}

\newcommand{\IL}{{\rm IL}}
\newcommand{\CL}{{\rm CL}}
\newcommand{\ILL}{{\rm AL}_i}
\newcommand{\CLL}{{\rm AL}_c}
\newcommand{\ILuL}{{\rm {\L}L}_i}
\newcommand{\CLuL}{{\rm {\L}L}_c}

\newcommand{\gbang}{{}!_\text{\rm g}}
\newcommand{\dbang}{{}!_\text{\rm d}}
\newcommand{\kbang}{{}!_\text{\rm k}}
\newcommand{\tkbang}{{}!_\text{\rm kt}}
\newcommand{\tdbang}{{}!_\text{\rm dt}}

\newcommand{\Aat}{A_{{\sf at}}}

\newcommand{\proves}{\vdash}
\newcommand{\eqleft}[1]{\begin{itemize} \item[] $#1$ \end{itemize}}

\newcommand{\uInter}[3]{|#1|^{#2}_{#3}}
\newcommand{\pvec}[1]{\boldsymbol{#1}}
\newcommand{\pdefin}{\mathrel{\mathop:}\equiv}
\newcommand{\setComp}[2]{\{ #1 \;\, | \;\, #2\}}

\newcommand{\kleene}[1]{\{ #1 \}}
\newcommand{\defined}[1]{{#1\! \downarrow}}

\newcommand{\NN}{\mathbb{N}}

\newcommand{\cwedge}{\otimes}
\newcommand{\scwedge}{\;\otimes\;}
\newcommand{\awedge}{\,\&\,}

\newcommand{\bang}{!}

\newcommand{\lNeg}{^\perp}
\newcommand{\lto}{\multimap}

\newcommand{\pcond}[3]{#2 \; \Diamond_{#1} \, #3}

\newcommand{\starT}[1]{#1^*}
\newcommand{\circT}[1]{#1^\circ}
\newcommand{\forgetT}[1]{#1^F}

\newcommand{\real}[2]{#2_{\sf r}(#1)}
\newcommand{\kreal}[2]{#2_{\sf nr}(#1)}
\newcommand{\rkreal}[3]{#3_{\sf rnr}(#1;#2)}
\newcommand{\rreal}[3]{#3_{\sf rr}(#1;#2)}
\newcommand{\qrreal}[3]{#3_{\sf qr}(#1;#2)}
\newcommand{\mrt}[2]{#2_{\sf mrt}(#1)}
\newcommand{\dial}[3]{#3_{\sf D}(#1;#2)}

\newcommand{\bQuant}[3]{\forall #1 \!\prec\! #2 \, #3}

\begin{document} 

%\title{A Few Lessons from Unifying Functional Interpretations}
\title{Unifying Functional Interpretations: Past and Future}

\author{Paulo Oliva \\[3mm]
Queen Mary University of London \\[2mm] \emph{paulo.oliva@eecs.qmul.ac.uk} \\[2mm] 
}

\maketitle
\thispagestyle{empty}

\begin{abstract} This article surveys work done in the last six years on the unification of various functional interpretations including G\"odel's dialectica interpretation, its Diller-Nahm variant, Kreisel modified realizability, Stein's family of functional interpretations, functional interpretations ``with truth", and bounded functional interpretations. Our goal in the present paper is twofold: (1) to look back and single out the main lessons learnt so far, and (2) to look forward and list several open questions and possible directions for further research.
\end{abstract} 

%%%%%%%%%%%%%%%%%%%%%%%%%%%%%%%%%%%%%%%%
%%%%%%%%%%%%%%%%%%%%%%%%%%%%%%%%%%%%%%%%
\section{Introduction}
%%%%%%%%%%%%%%%%%%%%%%%%%%%%%%%%%%%%%%%%
%%%%%%%%%%%%%%%%%%%%%%%%%%%%%%%%%%%%%%%%

When studying and working with the two main functional interpretations, namely the \emph{dialectica} \cite{Avigad(98),Goedel(58)} and the modified realizability \cite{Kreisel(59)} interpretations, one notices a striking similarity in the way the two interpretations behave. For instance, they both interpret $\forall \exists$-statements in precisely the same way, and their soundness (also called adequacy) proofs follow very similar patterns. Yet, for all purpose these are two very different interpretations, validating different principles\footnote{For instance, the dialectica interpretation validates the Markov principle whereas modified realizability does not. On the other hand, modified realizability validates full extensionality whereas the dialectica interpretation does not.}, and having different properties\footnote{For instance, realizability interpretations always have a so-called ``truth" variant, whereas the dialectica interpretation does not.}. Several questions naturally arise. What is the common structure behind these two functional interpretations? How are the different witnesses obtained from a given proof when applying different interpretations related to each other?

%	\item[({\bf Q2})] Is there is a common abstract interpretation that generalises both the dialectica and modified realizability?
%	\item[({\bf Q3})] 
%\end{itemize} 

It was with these questions in mind that I set out \cite{Oliva(06)} to develop a general framework to unify functional interpretations. This initial work was followed by several other articles \cite{FO(2009A),FO(2012),FO(2011A),GO(2010),Hernest(2008),Oliva(2007A),Oliva(2007),Oliva(2008),Oliva(2009A),Oliva(2009B)} further refining or generalising the original idea. These were mainly done in collaboration with Gilda Ferreira, Jaime Gaspar and Mircea-Dan Hernest. What started as a small modification of the dialectica interpretation to also capture realizability and the Diller-Nahm variant \cite{Diller(74)} ended up as a very general hybrid functional interpretation of intuitionistic affine logic\footnote{Intuitionistic linear logic plus the weakening rule.}, also capturing Stein's family of functional interpretations \cite{Stein(79)}, functional interpretations ``with truth" \cite{GO(2010)}, and bounded functional interpretations \cite{Ferreira(05A),Ferreira(05),Ferreira(06)}. 

This article will survey the work mentioned above, singling out what I believe to be the key lessons learnt so far. These are summarised as follows. For details see the corresponding sections and the articles mentioned.
\begin{enumerate}
	\item[(\S \ref{newReal})] Modified realizability can also alternatively be presented as a \emph{relation} between potential witnesses and challenges, in a way very similar to the way the dialectica interpretation is presented. This is originally observed in \cite{Oliva(06)} and is key to extending realizability to affine logic \cite{Oliva(2007)}.
	\item[(\S \ref{factoring})] Most functional interpretations of \emph{intuitionistic logic} can be factored via \emph{affine logic}. More interestingly, all functional interpretations considered, when extended to affine logic, coincide in the pure fragment, where modalities are absent. This factorisation allows us to clearly see that the only difference between most of the functional interpretations is in the treatment of contraction, which in affine logic is captured by $\bang A$. Although this was originally done in the setting of \emph{classical} affine logic \cite{Hernest(2008),Oliva(2007A),Oliva(2007),Oliva(2009A)}, it turned out that \emph{intuitionistic} affine logic is not only enough, but the unification becomes much simpler \cite{FO(2009A),FO(2011A),GO(2010)} (albeit at the cost of losing symmetry). 
	\item[(\S \ref{sec-truth})] When designing the unified functional interpretation of intuitionistic affine logic we were only expecting to be able to capture the classic interpretations such as the dialectica, modified realizability and Diller-Nahm. We were therefore surprised when we discovered \cite{GO(2010)} that even the truth variants of functional interpretations fit in the framework almost effortlessly. Which means that even proof interpretations with truth only differ from their ``non-truth" variants in the treatment of $\bang A$, but coincide in the treatment of all other connectives. 
	\item[(\S \ref{hybrid})] Because the bang ($\bang$) of affine logic is not canonical, one can then effectively combine all the functional interpretations mentioned above, including their truth variants, into single interpretations which we called \emph{hybrid functional interpretations} \cite{Hernest(2008),Oliva(2009B)}. This means, for instance, that in a single proof one can try to make use of both the dialectica interpretation in some parts of the proof and modified realizability in others, combining their strengths to maximum benefit. 
\end{enumerate}
We will conclude (\S \ref{future-work}) by listing thirteen open questions which indicate possible interesting directions for further research.

\vspace{2mm}

\noindent {\bf Acknowledgement}. Most of the work presented here has been done in collaboration with Gilda Ferreira, Jaime Gaspar and Mircea-Dan Hernest. I would also like to acknowledge previous work done in this direction on which the current work builds, such as those of Martin Stein \cite{Stein(79),Stein(1980)}, Val\'eria de Paiva \cite{dePaiva(1989A),dePaiva(1989B)}, Masaru Shirahata \cite{Shirahata(2006)} and Andreas Blass \cite{Blass(1992)}. Finally, many thanks to Thomas Powell, Jules Hedges and Gilda Ferreira for several comments and corrections on an earlier version of this paper.

\vspace{2mm}

\noindent {\bf Notation}. We use $X \pdefin A$ to say that $X$ is defined by $A$. We use $A \equiv B$ to mean $A$ and $B$ are syntactically equal.

%%%%%%%%%%%%%%%%%%%%%%%%%%%%%%%%%%%%%%%%
%%%%%%%%%%%%%%%%%%%%%%%%%%%%%%%%%%%%%%%%
\section{A Different View on Realizability}
%%%%%%%%%%%%%%%%%%%%%%%%%%%%%%%%%%%%%%%%
%%%%%%%%%%%%%%%%%%%%%%%%%%%%%%%%%%%%%%%%
\label{newReal}

The first obvious difference between modified realizability \cite{Kreisel(59)} and the dialectica interpretation \cite{Goedel(58)} is that the first interprets formulas $A$ as \emph{unary} predicates $\real{\pvec x}{A}$, normally written as ``$\pvec x$ realizes $A$", whereas the dialectica interpretation associates to formulas $A$ \emph{binary} predicates $\dial{\pvec x}{\pvec y}{A}$. Here $\pvec x$ and $\pvec y$ denote tuples of variables $\pvec x = x_1, \ldots, x_n$ and $\pvec y = y_1, \ldots, y_m$, where the length of the tuple and the types of the variables depend on the logical structure of the formula $A$. The two formulas $\real{\pvec x}{A}$ and $\dial{\pvec x}{\pvec y}{A}$ are defined inductively as\footnote{We are using the abbreviation $\pcond{b}{A}{B} \pdefin (b = {\sf true} \to A) \wedge (b = {\sf false} \to B)$. We also use the same macro in the context of affine logic where it stands for $\pcond{b}{A}{B} \pdefin (\bang (b = {\sf true}) \lto A) \cwedge (\bang (b = {\sf false}) \lto B)$. }
\[
{\small 
\begin{array}{lcl}
	\real{\pvec x, \pvec y}{(A \wedge B)} & \pdefin & \real{\pvec x}{A} \wedge \real{\pvec y}{B} \\[1mm]
	\real{\pvec x, \pvec y, b}{(A \vee B)} & \pdefin & \pcond{b}{\real{\pvec x}{A}}{\real{\pvec y}{B}} \\[1mm]
	\real{\pvec f}{(A \to B)} & \pdefin & \forall \pvec x (\real{\pvec x}{A} \to \real{\pvec f \pvec x}{B}) \\[1mm]
	\real{\pvec x, a}{(\exists z A)} & \pdefin & \real{\pvec x}{(A[a/z])} \\[1mm]
	\real{\pvec f}{(\forall z A)} & \pdefin & \forall z \real{\pvec f z}{A}
\end{array}
\quad
\begin{array}{lcl}
	\dial{\pvec x,\pvec v}{\pvec y, \pvec w}{(A \wedge B)} & \pdefin & \dial{\pvec x}{\pvec y}{A} \wedge \dial{\pvec v}{\pvec w}{B} \\[1mm]
	\dial{\pvec x,\pvec v,b}{\pvec y, \pvec w}{(A \vee B)} & \pdefin & \pcond{b}{\dial{\pvec x}{\pvec y}{A}}{\dial{\pvec v}{\pvec w}{B}} \\[1mm]
	\dial{\pvec f, \pvec g}{\pvec x, \pvec w}{(A \to B)} & \pdefin & \dial{\pvec x}{\pvec g \pvec x \pvec w}{A} \to \dial{\pvec f \pvec x}{\pvec w}{B} \\[1mm]
	\dial{\pvec x, a}{\pvec y}{(\exists z A)} & \pdefin & \dial{\pvec x}{\pvec y}{(A[a/z])} \\[1mm]
	\dial{\pvec f}{\pvec y, a}{(\forall z A)} & \pdefin & \dial{\pvec f a}{\pvec y}{(A[a/z])}.
\end{array}
}
\]
with the base case $\real{\epsilon}{(\Aat)} = \dial{\epsilon}{\epsilon}{(\Aat)} = \Aat$, for atomic formulas $\Aat$, with $\epsilon$ denoting the empty tuple (henceforth omitted). Note that for tuples of variables $\pvec f = f_1, \ldots, f_n$ and $\pvec x$ we write $\pvec f \pvec x$ for the tuple of terms $f_1 \pvec x, \ldots, f_n \pvec x$. Using these predicates $\real{\pvec x}{A}$ and $\dial{\pvec x}{\pvec y}{A}$ we can define two sets of ``functionals"
\[  A \quad \mapsto \quad \setComp{\pvec x}{\real{\pvec x}{A}} \hspace{2cm} A \quad \mapsto \quad \setComp{\pvec x}{\forall \pvec y \dial{\pvec x}{\pvec y}{A}} \]
which we will refer to as the ``realizability witnesses" and the ``dialectica witnesses". The two functional interpretations, modified realizability and dialectica, can be viewed as algorithms to turn an intuitionistic proof of $A$ into concrete (e.g. higher-order programs) elements of these sets.

The work on unifying different functional interpretations \cite{Oliva(06)} started with the observation that one can also view modified realizability as associating formulas with a \emph{binary} predicate $\rreal{\pvec x}{\pvec y}{A}$ (which I will call ``relational realizability") between two tuples $\pvec x$ and $\pvec y$ in a way very similar to the dialectica interpretation, namely
\begin{equation} \label{relational-realizability}
\begin{array}{lcl}
	\rreal{\pvec x, \pvec v}{\pvec y, \pvec w}{(A \wedge B)} & \pdefin & \rreal{\pvec x}{\pvec y}{A} \wedge \rreal{\pvec v}{\pvec w}{B} \\[1mm]
	\rreal{\pvec x, \pvec v, b}{\pvec y, \pvec w}{(A \vee B)} & \pdefin & \pcond{b}{\rreal{\pvec x}{\pvec y}{A}}{\rreal{\pvec v}{\pvec w}{B}} \\[1mm]
	\rreal{\pvec f}{\pvec x, \pvec w}{(A \to B)} & \pdefin & \forall \pvec y \rreal{\pvec x}{\pvec y}{A} \to \rreal{\pvec f \pvec x}{\pvec w}{B} \\[1mm]
	\rreal{\pvec x, a}{\pvec y}{(\exists z A)} & \pdefin & \rreal{\pvec x}{\pvec y}{(A[a/z])} \\[1mm]
	\rreal{\pvec f}{\pvec y, a}{(\forall z A)} & \pdefin & \rreal{\pvec f a}{\pvec y}{(A[a/z])}.
\end{array}
\end{equation}
It is easy to show by induction on the formula $A$ that these two different definitions of realizability lead to the same interpretation as the following equivalence is intuitionistically provable:
\begin{equation*}
\real{\pvec x}{A} \quad \Leftrightarrow \quad \forall \pvec y \rreal{\pvec x}{\pvec y}{A}.
\end{equation*}
The relational presentation of realizability, however, makes it absolutely clear that realizability only differs from the dialectica interpretation in the clause for implication $A \to B$. While the realizability interpretation does not attempt to witness the universal quantifier $\forall \pvec y$ in the clause for $A \to B$, the dialectica interpretation witnesses such quantifier via the extra tuple of functionals $\pvec g$. 

The two main ideas behind the original \emph{unifying functional interpretation} \cite{Oliva(06)} are the introduction of a common notation $\uInter{A}{\pvec x}{\pvec y}$ for such binary predicates, and a parametrised interpretation of $A \to B$. That is achieved via an abstract formula constructor $\bQuant{\pvec x}{\pvec a}{A}$ that takes a tuple of terms $\pvec a$ and a formula $A$ (with free variables $\pvec x$) and produces a new formula where $\pvec x$ are no longer free. A parametrised functional interpretation can then be given as
\begin{equation} \label{ufi-implication}
\begin{array}{lcl}
	\uInter{A \wedge B}{\pvec x, \pvec v}{\pvec y, \pvec w}& \pdefin &
		\uInter{A}{\pvec x}{\pvec y} \wedge \uInter{B}{\pvec v}{\pvec w} \\[1mm]
	\uInter{A \vee B}{\pvec x, \pvec v, b}{\pvec y, \pvec w} & \pdefin &
		\pcond{b}{\uInter{A}{\pvec x}{\pvec y}}{\uInter{B}{\pvec v}{\pvec w}} \\[1mm]
	\uInter{A \to B}{\pvec f, \pvec g}{\pvec x, \pvec w} & \pdefin &
		\bQuant{\pvec y}{\pvec g \pvec x \pvec w}{\uInter{A}{\pvec x}{\pvec y}} \to \uInter{B}{\pvec f \pvec x}{\pvec w} \\[1mm]
	\uInter{\exists z A}{\pvec x, a}{\pvec y} & \pdefin &
		\uInter{A[a/z]}{\pvec x}{\pvec y} \\[1mm]
	\uInter{\forall z A}{\pvec f}{\pvec y, a} & \pdefin &
		\uInter{A[a/z]}{\pvec f a}{\pvec y}.
\end{array}
\end{equation}
Subject to a few conditions (cf. \cite{Oliva(06)}) on $\bQuant{\pvec x}{\pvec a}{A}$, one can then prove a uniform soundness theorem for intuitionistic logic. When the formula constructor is instantiated one obtains the three main functional interpretations as follows:
\begin{equation*}
\begin{array}{c|c}
	\bQuant{\pvec x}{\pvec a}{A} & \mbox{{\bf Functional interpretation}} \\[1mm]
	\hline \\[-3mm]
	A[\pvec a/\pvec x] & \mbox{G\"odel's dialectica interpretation} \\
	\forall \pvec x \!\in\! \pvec a \, A & \mbox{Diller-Nahm interpretation} \\
	\forall \pvec x A & \mbox{Kreisel modified realizability}
\end{array}
\end{equation*}
In order to show that each of these three interpretations is sound one only needs to check that they satisfy the required conditions mentioned above.

\begin{remark}[Stein family of interpretations] Let $M \in \NN \cup \{ \infty \}$. Given a tuple of variables $\pvec x = x_0, \ldots, x_n$ let us denote by $\pvec x^{\geq M}$ the tuple containing only the elements of $\pvec x$ with type level $\geq M$. Similarly we denote by $\pvec x^{< M}$ the tuple containing only the elements of $\pvec x$ with type level $< M$. Note that $\pvec x^{< \infty} = \pvec x$ and $\pvec x^{< 0}$ is the empty tuple. Stein's family of functional interpretations \cite{Stein(79)} also fits in the above framework as we can take for each given $M$
\[ \bQuant{\pvec x}{\pvec a}{A} \pdefin \forall \pvec x^{< M} \forall \pvec x^{\geq M} \!\in\! \pvec a \, A \]
where $\pvec a$ is a set indexed by the pure type $M$, i.e. $\pvec a \colon M \to \rho$ for some type $\rho$. When $M = \infty$ this coincides with modified realizability, whereas with $M = 0$ this is a variant of the Diller-Nahm interpretation that allows for infinite (countable) sets, as $\pvec a \colon \NN \to \rho$ ($\NN$ is the pure type having type level $0$).
%This observation led Stein \cite{Stein(1980)} to develop the first unification of functional interpretation based a number parameter $M$ determining the type level of variables in the tuple $\pvec y$ which should and should not be witnessed. Stein's idea is to interpret implication and universal quantifiers as
%%
%\eqleft{
%\begin{array}{lcl}
%\stein{\pvec f, \pvec g}{\pvec x^{\geq M}, \pvec w}{(A \to B)} & \pdefin & \forall \pvec x^{< M} ( \forall \pvec y \in \pvec g \pvec x \pvec w \, \stein{\pvec x}{\pvec y}{A} \to \stein{\pvec f \pvec x}{\pvec w}{B}), \\[2mm]
%\stein{\pvec f}{\pvec y^{\geq M}, a}{(\forall z A)} & \pdefin & \forall \pvec y^{< M} \stein{\pvec f a}{\pvec y}{(A[a/z])}.
%\end{array}
%}
%%
%The bounded quantification $\forall \pvec y \!\in\! \pvec a \, A$ is given a direct axiomatisation sufficient for the soundness interpretation to go through (rather than being defined via a primitive relation $\pvec y \in \pvec a$). 
%
%For technical reasons Stein does not treat G\"odel's original dialectica interpretation, but rather its Diller-Nahm variant \cite{Diller(74)}. The general unification proposed in \cite{Oliva(06)} extends Stein's unification and easily covers the dialectica interpretation as well\footnote{At the time of writing \cite{Oliva(06)} I was unaware of Stein's similar unification of functional interpretations in \cite{Stein(1980)}. I would like to thank Jaime Gaspar for bringing \cite{Stein(1980)} to my attention in June 2010.}.
\end{remark}

%%%%%%%%%%%%%%%%%%%%%%%%%%%%%%%%%%%%%%%%
%%%%%%%%%%%%%%%%%%%%%%%%%%%%%%%%%%%%%%%%
\section{Factoring Through Affine Logic}
%%%%%%%%%%%%%%%%%%%%%%%%%%%%%%%%%%%%%%%%
%%%%%%%%%%%%%%%%%%%%%%%%%%%%%%%%%%%%%%%%
\label{factoring}

Reformulating realizability as a binary predicate as described in Section \ref{newReal} was an important step towards showing that modified realizability and the dialectica interpretation have much more in common than previously imagined. The fact is that they only differ on their handling of witnesses coming from the premise of an implication. But that opens a new question: What is special about the premise of an implication that allows for these different interpretations to exist? A satisfactory answer to this question came from the analysis of functional interpretations via affine logic.

Intuitionistic affine logic ($\ILL$) is a refinement of intuitionistic logic ($\IL$) where particular attention is paid to the contraction rule \cite{Benton(1993),Girard(87B)}. We call this a \emph{refinement} because the connectives of intuitionistic logic can be recovered from a combination of those from affine logic. This is formally expressed via Girard's translations of intuitionistic logic into linear logic. The two most commonly used are\footnote{The usual clause for $\starT{(A \wedge B)}$ is $\starT{(A \wedge B)} \pdefin \starT{A} \awedge \starT{B}$. We can take $\starT{(A \wedge B)} \pdefin \starT{A} \cwedge \starT{B}$ instead because we are embedding intuitionistic logic into affine logic (linear logic with the weakening rule).}
\[
\begin{array}{lcl}
\starT{P} & \pdefin & P \\[2mm]
\starT{(A \wedge B)} & \pdefin & \starT{A} \cwedge \starT{B} \\[2mm]
\starT{(A \vee B)} & \pdefin & \bang \starT{A} \, \oplus \; \bang \starT{B} \\[2mm]
\starT{(A \to B)} & \pdefin & \bang \starT{A} \lto \starT{B} \\[2mm]
\starT{(\forall x A)} & \pdefin & \forall x \starT{A} \\[2mm]
\starT{(\exists x A)} & \pdefin & \exists x \bang \starT{A}.
\end{array}
\hspace{10mm}
\begin{array}{lcl}
\circT{P} & \pdefin & \bang P \\[2mm]
\circT{(A \wedge B)} & \pdefin & \circT{A}\cwedge \circT{B} \\[2mm]
\circT{(A \vee B)} & \pdefin & \circT{A} \oplus \circT{B} \\[2mm]
\circT{(A \to B)} & \pdefin & \bang (\circT{A} \lto \circT{B}) \\[2mm]
\circT{(\forall x A)} & \pdefin & \bang \forall x \circT{A} \\[2mm]
\circT{(\exists x A)} & \pdefin & \exists x \circT{A}.
\end{array}
\]
The translations are such that if $A$ is provable in $\IL$ then both $\bang \starT{A}$ and $\circT{A}$ are provable in $\ILL$. 
% Clearly one can always turn a proof in linear logic into a proof in classical logic by replacing $\cvee, \oplus$ by disjunction, $\cwedge, \awedge$ by conjunction, and omitting the modalities. We denote this translation of LL as $A^F$.

While working on \cite{Oliva(06)}, in the setting of intuitionistic logic, I came across de Paiva's \cite{dePaiva(1989B)} dialectica (and Diller-Nahm) interpretation of \emph{affine logic}. It then occurred to me that one could use the new formulation of realizability discussed in Section \ref{newReal} to extend the \emph{realizabillity} interpretation from intuitionistic logic to affine logic. This was developed and presented in \cite{Oliva(2007A),Oliva(2007)}. The starting point is the functional interpretation of pure affine logic (affine logic without the exponentials). As mentioned in the introduction, we consider the intuitionistic fragment of affine logic:
\begin{equation} \label{basic-inter}
\begin{array}{lcl}
%
%\uInter{A \awedge B}{\pvec x, \pvec v}{\pvec y, \pvec w,z} & \pdefin & \pcond{z}{\uInter{A}{\pvec x}{\pvec y}}{\uInter{B}{\pvec v}{\pvec w}} \\[2mm]
%
\uInter{A \oplus B}{\pvec x, \pvec v,z}{\pvec y, \pvec w} & \pdefin & \pcond{z}{\uInter{A}{\pvec x}{\pvec y}}{\uInter{B}{\pvec v}{\pvec w}} \\[2mm]
\uInter{A \cwedge B}{\pvec x, \pvec v}{\pvec y, \pvec w} & \pdefin & \uInter{A}{\pvec x}{\pvec y} \cwedge \uInter{B}{\pvec v}{\pvec w} \\[2mm]
\uInter{A \lto B}{\pvec f, \pvec g}{\pvec x, \pvec w} & \pdefin & \uInter{A}{\pvec x}{\pvec g \pvec x \pvec w} \lto \uInter{B}{\pvec f \pvec y}{\pvec w} \\[2mm]
\uInter{\forall z A(z)}{\pvec f}{\pvec y, a} & \pdefin & \uInter{A[a/z]}{\pvec f a}{\pvec y} \\[2mm]
\uInter{\exists z A(z)}{\pvec x, a}{\pvec y} & \pdefin & \uInter{A[a/z]}{\pvec x}{\pvec y}.
\end{array}
\end{equation}

\begin{figure}[t]
\begin{center}
	\setlength{\unitlength}{10mm}
	\begin{picture}(6.0,3.3)
		% border
		% \put(0.0,0.0){\line(1,0){6}}
		% \put(0.0,3.0){\line(1,0){6}}
		\thicklines
		\put(0.5,2.8){$A$}
		\put(0.65,2.65){\vector(0,-1){1.8}}
			\put(-0.2,1.7){$\circT{(\cdot)}$}
		\put(1.1,2.95){\vector(1,0){3.3}}
			\put(1.5,2.5){mod. realizability}
		\put(0.5,0.3){$\circT{A}$}
		\put(4.5,2.8){$\real{\pvec x}{A}$}
		\put(4.65,2.6){\vector(0,-1){1.8}}
			\put(5.1,1.7){$\circT{(\cdot)}$}
		\put(1.2,0.45){\vector(1,0){3.2}}
			\put(2.4,0.7){$\uInter{\cdot}{}{}$}
			% \put(1.5,-0.3){$\bQuant{\pvec x}{\pvec a}{A} \equiv \forall \pvec x A$}
			\put(2.0,0.0){(\ref{basic-inter}) + (\ref{bang-realizability})}
		\put(4.5,0.3){$\uInter{\circT{A}}{\pvec x}{} \equiv \circT{(\real{\pvec x}{A})}$}
	\end{picture}
\caption{Factoring modified realizability}
\label{factor-realizability}
\end{center}
\end{figure}

What one notices is that the parameter constructor $\bQuant{\pvec x}{\pvec a}{A}$ used to interpret $A \to B$ in (\ref{ufi-implication}) is in fact the interpretation of the affine logic modality $\bang A$.
So we can extend the basic interpretation (\ref{basic-inter}) to a parametrised interpretation of full intuitionistic affine logic as
\begin{equation} \label{bang-unifying}
\begin{array}{lcl}
\uInter{\bang A}{\pvec x}{\pvec a} & \pdefin & \bang \bQuant{\pvec y}{\pvec a}{\uInter{A}{\pvec x}{\pvec y}}.
\end{array}
\end{equation}
Via the translations $\circT{(\cdot)}$ and $\starT{(\cdot)}$ of $\IL$ into $\ILL$ one can recover the interpretations of intuitionistic logic from those of intuitionistic affine logic as follows. For instance, consider the abbreviation $\bQuant{\pvec x}{\pvec a}{A} \pdefin \forall \pvec x A$, so that (\ref{bang-unifying}) simplifies to
\begin{equation} \label{bang-realizability}
\begin{array}{lcl}
\uInter{\bang A}{\pvec x}{} & \pdefin & \bang \forall \pvec y \uInter{A}{\pvec x}{\pvec y}.
\end{array}
\end{equation}
We call the resulting interpretation a \emph{modified realizability} interpretation of affine logic because the diagram of Figure \ref{factor-realizability} commutes, i.e. given a formula $A$ of intuitionistic logic we can either apply modified realizability directly and translate the result into liner logic, or alternatively, we can first translate $A$ into affine logic, and then apply the interpretation with $\bQuant{\pvec x}{\pvec a}{A} \pdefin \forall \pvec x A$. Both paths result in the \emph{same} formula. Note that we really mean syntactic equality, rather than logical equivalence.

Now, if instead of using the Girard translation $\circT{A}$ we use instead the translation $\starT{A}$ we obtain a different diagram (Figure \ref{factor-relational-realizability}) which also commutes if we take in the upper arrow the \emph{relational} realizability instead.

\begin{figure}[h]
\begin{center}
	\setlength{\unitlength}{10mm}
	\begin{picture}(6.0,3.3)
		% border
		% \put(0.0,0.0){\line(1,0){6}}
		% \put(0.0,3.0){\line(1,0){6}}
		\thicklines
		\put(0.5,2.8){$A$}
		\put(0.65,2.65){\vector(0,-1){1.8}}
			\put(-0.2,1.7){$\starT{(\cdot)}$}
		\put(1.1,2.95){\vector(1,0){3.3}}
			\put(1.5,2.5){rel. realizability}
		\put(0.5,0.3){$\starT{A}$}
		\put(4.5,2.8){$\rreal{\pvec x}{\pvec y}{A}$}
		\put(4.65,2.65){\vector(0,-1){1.8}}
			\put(5.1,1.7){$\starT{(\cdot)}$  \quad (cf. Remark \ref{remark-trans-simp} below)}
		\put(1.2,0.45){\vector(1,0){3.2}}
			\put(2.4,0.7){$\uInter{\cdot}{}{}$}
			% \put(1.5,-0.3){$\bQuant{\pvec x}{\pvec a}{A} \pdefin \forall \pvec x A$}
			\put(2.0,0.0){(\ref{basic-inter}) + (\ref{bang-realizability})}
		\put(4.5,0.3){$\uInter{\starT{A}}{\pvec x}{\pvec y} \equiv \starT{(\rreal{\pvec x}{\pvec y}{A})}$}
	\end{picture}
\caption{Factoring the relational variant of modified realizability}
\label{factor-relational-realizability}
\end{center}
\end{figure}

In other words, the two ways of presenting modified realizability arise from the two possible ways of translating intuitionistic logic into affine logic. In both cases the modified realizability interpretation of \emph{affine logic} is fixed (the lower arrows of Figures \ref{factor-realizability} and \ref{factor-relational-realizability}). That illustrates how affine logic has a more fundamental nature, as it is able to capture precisely the inherent structure of realizability. 

Just as we have factored the realizability interpretation through affine logic, we can also do the same for the dialectica interpretation by considering the abbreviation $\bQuant{\pvec x}{\pvec a}{A} \pdefin A[\pvec a/\pvec x]$ leading to the interpretation of $\bang A$ as 
\begin{equation} \label{bang-dialectica}
\begin{array}{lcl}
\uInter{\bang A}{\pvec x}{\pvec y} & \pdefin & \bang \uInter{A}{\pvec x}{\pvec y}.
\end{array}
\end{equation}
Again, we say that (\ref{bang-dialectica}) is a dialectica interpretation of affine logic because it corresponds to the dialectica interpretation of intuitionistic logic as depicted in the commuting diagram of Figure \ref{factor-dialectica}.

\begin{figure}[h]
\begin{center}
	\setlength{\unitlength}{10mm}
	\begin{picture}(6.0,3.3)
		% border
		% \put(0.0,0.0){\line(1,0){6}}
		% \put(0.0,3.0){\line(1,0){6}}
		\thicklines
		\put(0.5,2.8){$A$}
		\put(0.65,2.65){\vector(0,-1){1.8}}
			\put(-0.2,1.7){$\starT{(\cdot)}$}
		\put(1.1,2.95){\vector(1,0){3.3}}
			\put(1.9,2.5){dialectica}
		\put(0.5,0.3){$\starT{A}$}
		\put(4.5,2.8){$\dial{\pvec x}{\pvec y}{A}$}
		\put(4.65,2.65){\vector(0,-1){1.8}}
			\put(5.1,1.7){$\starT{(\cdot)}$ \quad (cf. Remark \ref{remark-trans-simp} below)}
		\put(1.2,0.45){\vector(1,0){3.2}}
			\put(2.4,0.7){$\uInter{\cdot}{}{}$}
			% \put(1.4,-0.3){$\bQuant{x}{a}{A} \pdefin A[a/x]$}
			\put(2.0,0.0){(\ref{basic-inter}) + (\ref{bang-dialectica})}
		\put(4.5,0.3){$\uInter{\starT{A}}{\pvec x}{\pvec y} \equiv \starT{(\dial{\pvec x}{\pvec y}{A})}$}
	\end{picture}
\caption{Factoring G\"odel's dialectica interpretation}
\label{factor-dialectica}
\end{center}
\end{figure}

Finally, a Diller-Nahm interpretation of affine logic is obtained by choosing the abbreviation
\[ \bQuant{\pvec x}{\pvec a}{A} \pdefin \forall \pvec x \!\in\! \pvec a \, A, \]
where $\pvec a$ is a tuple of finite sets, and $\pvec x \in \pvec a$ denotes the usual set inclusion. For further details on the factorisation of the main functional interpretations via affine logic see \cite{FO(2009A),FO(2011A)}.

\begin{remark} \label{remark-trans-simp} In the diagrams of Figures \ref{factor-relational-realizability} and \ref{factor-dialectica} we are taking a simplified form of the $\starT{(\cdot)}$-translation, namely, one where the clauses for disjunction and existential quantifier are simply
\[
\begin{array}{lcl}
\starT{(A \vee B)} & \pdefin & \starT{A} \oplus \starT{B} \\[2mm]
\starT{(\exists x A)} & \pdefin & \exists x \starT{A},
\end{array}
\]
i.e. the bang is not used. The reason why we can work with this simpler translation of $\IL$ into $\ILL$ is because we are considering $\ILL$ extended with the following two principles
\begin{equation} \label{ufi-extra-principles}
\begin{array}{c}
\bang A \, \oplus \, \bang B \lto~\bang (A \oplus B) \\[2mm]
\exists x \bang A \lto~\bang \exists x A.
\end{array}
\end{equation}
These principles are harmless because they are interpretable by the interpretation $\uInter{A}{\pvec x}{\pvec y}$ for any of the three choices of $\bQuant{\pvec x}{\pvec a}{A}$ above. In general however, the combination of $\uInter{A}{\pvec x}{\pvec y}$ with the translation $\starT{(\cdot)}$ will lead to interpretations of disjunction and existential quantifier as
\begin{equation} \label{ufi-exists}
\begin{array}{lcl}
\uInter{A \vee B}{\pvec x, \pvec v, b}{\pvec y, \pvec w} & \pdefin & \pcond{b}{\bQuant{\pvec y}{\pvec a}{\uInter{A}{\pvec x}{\pvec y}}}{\bQuant{\pvec w}{\pvec c}{\uInter{B}{\pvec v}{\pvec w}}} \\[2mm]
\uInter{\exists z A}{\pvec x, a}{\pvec c} & \pdefin & \bQuant{\pvec y}{\pvec c}{\uInter{A[a/z]}{\pvec x}{\pvec y}}.
\end{array}
\end{equation}
This more general treatment is important for instance in the functional interpretation with truth as discussed in the following section.
\end{remark}

%%%%%%%%%%%%%%%%%%%%%%%%%%%%%%%%%%%%%%%%
%%%%%%%%%%%%%%%%%%%%%%%%%%%%%%%%%%%%%%%%
\section{Interpretations with Truth}
%%%%%%%%%%%%%%%%%%%%%%%%%%%%%%%%%%%%%%%%
%%%%%%%%%%%%%%%%%%%%%%%%%%%%%%%%%%%%%%%%
\label{sec-truth}

\begin{figure}[t]
\begin{center}
	\setlength{\unitlength}{10mm}
	\begin{picture}(6.0,3.3)
		% border
		% \put(0.0,0.0){\line(1,0){6}}
		% \put(0.0,3.0){\line(1,0){6}}
		\thicklines
		\put(0.0,2.8){$A$}
		\put(0.15,2.65){\vector(0,-1){1.8}}
			\put(-0.7,1.7){$\circT{(\cdot)}$}
		\put(0.6,2.95){\vector(1,0){4.3}}
			\put(1.0,2.5){realizability with truth}
		\put(0.0,0.3){$\circT{A}$}
		\put(5.0,2.8){$\mrt{\pvec x}{A}$}
		\put(5.15,2.65){\vector(0,-1){1.8}}
			\put(5.45,1.7){$\circT{(\cdot)}$}
		\put(0.7,0.45){\vector(1,0){4.2}}
			% \put(2.4,0.4){}
			\put(2.0,0.0){(\ref{basic-inter}) + (\ref{bang-real-truth})}
		\put(5.0,0.3){$\uInter{\circT{A}}{\pvec x}{} \equiv \circT{(\mrt{\pvec x}{A})}$}
	\end{picture}
\caption{Factoring modified realizability with truth}
\label{factor-realizability-with-truth}
\end{center}
\end{figure}

The soundness of functional interpretations guarantees that from a proof of $A$ a tuple of terms $\pvec t$ can be extracted such that $\uInter{A}{\pvec t}{\pvec y}$. An important issue is that such a tuple $\pvec t$ provides a witness to the statement $\exists \pvec x \forall \pvec y \uInter{A}{\pvec x}{\pvec y}$, but not necessarily a witness to the original theorem $A$. For realizability interpretations some variants have been developed so that a realiser for $\exists z A$ also contains a witness for $z$. These are the so-called q-\emph{realizability} and \emph{realizability with truth} \cite{Grayson(1981),Kleene(45),Troelstra(98)}. In general what we would like is that
\begin{equation} \label{truth-property}
\forall \pvec y \uInter{A}{\pvec x}{\pvec y} \to A
\end{equation}
is derivable without the need for the characterisation principles\footnote{The characterisation principles are the extra logical principles needed to show the equivalence between $A$ and its interpretation $\exists \pvec x \forall \pvec y \uInter{A}{\pvec x}{\pvec y}$.} of the interpretation $\uInter{\cdot}{}{}$, because then we can extract actual witnesses from proofs of existential statements as follows
\[ \proves \exists z A(z) \qquad \stackrel{\textup{soundness}}{\Rightarrow} \qquad \proves \uInter{\exists z A(z)}{\pvec t, s}{\pvec y} \qquad \stackrel{(\ref{basic-inter})}{\equiv} \qquad \proves \uInter{A(s)}{\pvec t}{\pvec y} \qquad \stackrel{(\ref{truth-property})}{\Rightarrow} \qquad \proves A(s). \]
In joint work with Jaime Gaspar \cite{GO(2010)} we have shown how interpretations with truth arise from a slight modification of the abstract interpretation of $\bang A$ from (\ref{bang-unifying}) to
\begin{equation} \label{bang-unifying-truth}
\begin{array}{lcl}
\uInter{\bang A}{\pvec x}{\pvec a} & \pdefin & \bang \bQuant{\pvec y}{\pvec a}{\uInter{A}{\pvec x}{\pvec y}} \;\cwedge\; \bang A.
\end{array}
\end{equation}
For instance, if we take the realizability abbreviation $\bQuant{\pvec y}{\pvec a}{A} \pdefin \forall \pvec y A$ in this case we obtain
\begin{equation} \label{bang-real-truth}
\uInter{\bang A}{\pvec x}{} \pdefin~\bang \forall \pvec y \uInter{A}{\pvec x}{\pvec y} \;\cwedge\; \bang A.
\end{equation}
The composition of this affine logic interpretation with the translation $\circT{(\cdot)}$ gives us precisely the \emph{modified realizability with truth} \cite{Kohlenbach(92A),Kohlenbach(98A),Kohlenbach(2008)}, as described in the diagram of Figure \ref{factor-realizability-with-truth}.

Consider then the q-variant of the relational realizability (\ref{relational-realizability}) where the clauses for disjunction and existential quantification are modified as
\begin{equation} \label{q-relational-realizability}
\begin{array}{lcl}
	\qrreal{\pvec x, \pvec v, b}{}{(A \vee B)} & \pdefin & \pcond{b}{(\forall \pvec y \qrreal{\pvec x}{\pvec y}{A} \wedge A)}{(\forall \pvec w \qrreal{\pvec v}{\pvec w}{B} \wedge B)} \\[2mm]
	\qrreal{\pvec x, a}{}{(\exists z A)} & \pdefin & \forall \pvec y \qrreal{\pvec x}{\pvec y}{(A[a/z])} \wedge A[a/z].
\end{array}
\end{equation}
The diagram of Figure \ref{factor-q-realizability} shows how such q-realizability corresponds to the $\starT{(\cdot)}$ translation, making use in this particular case of the forgetful translation $\forgetT{(\cdot)}$ of affine logic back into intuitionistic logic instead\footnote{In this case a diagram similar to the ones considered before would not lead to a commuting diagram (not even if logical equivalence is taken instead of syntactic equality). The problem is that whereas $A$ might contain existential quantifiers its interpretation $\qrreal{\pvec x}{\pvec y}{A}$ does not. Hence, formulas which are duplicated in $\uInter{\starT{A}}{\pvec x}{\pvec y}$ because of the $\bang$ in $\exists x \bang A$ are not duplicated in $\starT{(\qrreal{\pvec x}{\pvec y}{A})}$ because the existential quantifiers have disappeared. One way to solve this is presented in \cite{GO(2010)}, but uses logical equivalence. Here we present an alternative solution which is to use the forgetful translation that leads to a commuting diagram with syntactic equality instead. Obviously this is a weaker result than the previous four diagrams, as $\starT{(A^I)} \equiv (\starT{A})^J$ implies $A^I \equiv \forgetT{((\starT{A})^J)}$ but not conversely.}.

\begin{figure}[h]
\begin{center}
	\setlength{\unitlength}{10mm}
	\begin{picture}(6.0,3.3)
		% border
		% \put(0.0,0.0){\line(1,0){6}}
		% \put(0.0,3.0){\line(1,0){6}}
		\thicklines
		\put(0.0,2.8){$A$}
		\put(0.15,2.65){\vector(0,-1){1.8}}
			\put(-0.7,1.7){$\starT{(\cdot)}$}
		\put(0.6,2.95){\vector(1,0){4.3}}
			\put(1.5,2.5){q-realizability}
		\put(0.0,0.3){$\starT{A}$}
		\put(5.0,2.8){$\qrreal{\pvec x}{\pvec y}{A} \equiv \forgetT{(\uInter{\starT{A}}{\pvec x}{\pvec y})}$}
		\put(5.15,0.8){\vector(0,1){1.8}}
			\put(5.45,1.7){$\forgetT{(\cdot)}$}
		\put(0.7,0.45){\vector(1,0){4.2}}
			% \put(2.4,0.4){}
			\put(2.0,0.0){(\ref{basic-inter}) + (\ref{bang-real-truth})}
		\put(5.0,0.3){$\uInter{\starT{A}}{\pvec x}{\pvec y}$}
	\end{picture}
\caption{Factoring q-realizability}
\label{factor-q-realizability}
\end{center}
\end{figure}

\noindent If one observes that the $\circT{A}$ translation is affine logic equivalent to the ``banged" $\starT{A}$ translation, i.e. $\circT{A} \leftrightarrow~\bang \starT{A}$, one obtains the following interesting (apparently unobserved) correspondence between realizability with truth and q-realizability 
\[ \mrt{\pvec x}{A} \quad \stackrel{\IL}{\Leftrightarrow} \quad \forall \pvec y \qrreal{\pvec x}{\pvec y}{A} \wedge A. \]

A great benefit of this analysis of truth interpretations via affine logic is that it gave us a handle to define truth variants of other functional interpretations. For instance, contrary to what was thought \cite{Joergensen(2004)}, we can immediately obtain a Diller-Nahm with truth instantiating (\ref{bang-unifying-truth}) as 
\[ 
\uInter{\bang A}{\pvec x}{\pvec a} \pdefin~\bang \forall \pvec y \!\in\! \pvec a \, \uInter{A}{\pvec x}{\pvec y} \;\cwedge\; \bang A.
\]
For more details on the unification of functional interpretations with truth see  \cite{GO(2010)}.

\begin{remark} It is essential here that one uses the full $\starT{(\cdot)}$ translation, not the simplification of the previous section (cf. Remark \ref{remark-trans-simp}) as the choice of interpretation (\ref{bang-unifying-truth}) for $\bang A$, although sound for affine logic, it is not sound for the extra principles (\ref{ufi-extra-principles}). 
\end{remark}

%%%%%%%%%%%%%%%%%%%%%%%%%%%%%%%%%%%%%%%%
%%%%%%%%%%%%%%%%%%%%%%%%%%%%%%%%%%%%%%%%
\section{Putting it All Together}
%%%%%%%%%%%%%%%%%%%%%%%%%%%%%%%%%%%%%%%%
%%%%%%%%%%%%%%%%%%%%%%%%%%%%%%%%%%%%%%%%
\label{hybrid}

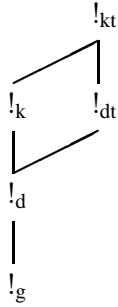
\begin{figure}[h]
\begin{center}
	\setlength{\unitlength}{8mm}
	\begin{picture}(4.5,4.5)
		\thicklines
		\put(1.0,-0.5){$\gbang$}
			\put(1.1,0.0){\line(0,1){0.7}}
		\put(1.0,1.0){$\dbang$}
			\put(1.1,1.5){\line(0,1){0.7}}
			\put(1.1,1.5){\line(2,1){1.4}}
		\put(1.0,2.5){$\kbang$}
			\put(1.1,3.0){\line(2,1){1.4}}
			\put(2.5,3.0){\line(0,1){0.7}}
		\put(2.4,2.5){$\tdbang$}
		\put(2.4,4.0){$\tkbang$}
	\end{picture}
\end{center}
\label{bang-diagram}
\caption{Ordering between different interpretations of $\bang A$}
\end{figure}

The analysis of different functional interpretations via \emph{affine logic} not only provides a setting where the precise differences between the interpretations can be clearly seen, but surprisingly it also allows us to combine multiple interpretations when analysing a single proof. This follows because, as observed by Girard, the bang ($\bang A$) is not a canonical operator. One can add multiple instances $\bang' A, \bang'' A, \ldots$ all with the same four rules without being able to show that any two are provably equivalent. This observation led us \cite{Hernest(2008),Oliva(2009B)} to consider a system of multi-modal affine logic with a different instance of $\bang A$ for each of the functional interpretations discussed above. For instance, we could add five different variants of $\bang A$ and interpret each as follows:
\[
\begin{array}{lll}
  \uInter{\kbang A}{\pvec x}{} & \pdefin~\bang \forall \pvec y \uInter{A}{\pvec x}{\pvec y} & \text{(Kreisel's modified realizability)} \\[2mm]
  \uInter{\dbang A}{\pvec x}{\pvec a} & \pdefin~\bang \forall \pvec y \!\in\! \pvec a \uInter{A}{\pvec x}{\pvec y} & \text{(Diller-Nahm interpretation)} \\[2mm]
  \uInter{\gbang A}{\pvec x}{\pvec y} & \pdefin~\bang \uInter{A}{\pvec x}{\pvec y} & \text{(G\"odel's dialectica interpretation)} \\[2mm]
  \uInter{\tkbang A}{\pvec x}{} & \pdefin~\bang \forall \pvec y \uInter{A}{\pvec x}{\pvec y} \cwedge~\bang A & \text{(Kreisel's modified realizability with truth)} \\[2mm]
  \uInter{\tdbang A}{\pvec x}{\pvec a} & \pdefin~\bang \forall \pvec y \!\in\! \pvec a \uInter{A}{\pvec x}{\pvec y} \cwedge~\bang A \qquad & \text{(Diller-Nahm interpretation with truth)}
\end{array}
\]
This leads to what we have termed \emph{hybrid functional interpretations}. If left completely unrelated, however, it would be difficult to make any practical use of this idea. We can observe, however, that there is a certain partial order between these different modalities, as for instance, a witness for $\kbang A$ is clearly also a witness for $\dbang A$. Therefore, we can add a rule that allows us to conclude $\dbang A$ from $\kbang A$, i.e.
\[
\begin{prooftree}
	\Gamma \proves \kbang A
	\justifies
	\Gamma \proves \dbang A 
\end{prooftree}
\]
In the diagram of Figure \ref{bang-diagram} we write $\bang_X$ above $\bang_Y$ if the interpretation of $\bang_X A$ implies the interpretation of $\bang_Y A$. As such, we could say that modified realizability with truth and G\"odel's \emph{dialectica} interpretation are the two ``extreme" interpretations amongst these five. For more details on these hybrid functional interpretations see \cite{GO(2010),Hernest(2008),Oliva(2009B)}.

%%%%%%%%%%%%%%%%%%%%%%%%%%%%%%%%%%%%%%%%
%%%%%%%%%%%%%%%%%%%%%%%%%%%%%%%%%%%%%%%%
\section{Directions for Further Work}
%%%%%%%%%%%%%%%%%%%%%%%%%%%%%%%%%%%%%%%%
%%%%%%%%%%%%%%%%%%%%%%%%%%%%%%%%%%%%%%%%
\label{future-work}

Let us conclude by outlining a few possible directions for further work. These are either directly related to the unification of functional interpretation or to the actual nature and better understanding of functional interpretations themselves. 

%{\bf Issues to be resolved}:
%
%\begin{itemize}
%	\item Why are there two ways of interpreting $\ILL$? Need to explain why in $\ILL$ one can either interpret implication and bang as
%	%
%	\[
%	\begin{array}{lcl}
%	\uInter{A \lto B}{\pvec f, \pvec g}{\pvec x, \pvec w} & \pdefin & \uInter{A}{\pvec x}{\pvec f \pvec x \pvec w} \lto \uInter{B}{\pvec g \pvec x}{\pvec w} \\[2mm]
%	\uInter{\bang A}{\pvec x}{\pvec a} & \pdefin & \forall \pvec y \in \pvec a \, \uInter{A}{\pvec x}{\pvec y}
%	\end{array}
%	\]
%	%
%	or in a more symmetric way as
%	%
%	\[
%	\begin{array}{lcl}
%	\uInter{A \lto B}{\pvec f, \pvec g}{\pvec x, \pvec w} & \pdefin & \uInter{A}{\pvec x}{\pvec f \pvec w} \lto \uInter{B}{\pvec g \pvec x}{\pvec w} \\[2mm]
%	\uInter{\bang A}{\pvec x}{\pvec f} & \pdefin & \forall \pvec y \in \pvec f \pvec x \, \uInter{A}{\pvec x}{\pvec y}.
%	\end{array}
%	\]
%	%
%	There is also the issue that $A \cwedge B$ has a complicated classical interpretation, but intuitionistically is very simple.
%	
%	\item What about completeness, can a nice pca be found where completeness holds, even for some fragment of the logic?
%	
%	\item What sets realizability apart?
%		
%	\item What has been learnt about interpretations with truth? Any other interpretations arise from general abstract view?
%
%\end{itemize}

%%%%%%%%%%%%%%%%%%%%%%%%%%%%%%%%%%%%%%%%
\subsection{Functional interpretations with forcing}
%%%%%%%%%%%%%%%%%%%%%%%%%%%%%%%%%%%%%%%%

The combination of realizability with Cohen's notion of forcing was originally studied by Goodman \cite{Goodman(1978)} who showed it to be an effective way to prove conservation results that cannot apparently be shown by realizability alone. Goodman's work is related to the interpretations with truth (cf. Section \ref{sec-truth}) as forcing is used precisely to recover the truth property (\ref{truth-property}). Although Goodman presented a single combined interpretation, Beeson \cite{Beeson(1979)} showed that Goodman's interpretation can actually be seen as a simple composition of the Kleene number realizability based on Turing machines with oracles followed by an application of forcing. Recently, another variant of realizability, called \emph{learning-based realizability} \cite{Aschieri(2010)}, has been developed providing an extension of realizability to classical arithmetic. Although different from Goodman's, the learning-based realizability has many similar features to Goodman's combination of realizability and forcing. For instance, the learning-based interpretation of formulas is described relative to a \emph{memory}, which can be understood as a forcing condition approximating a non-computable oracle. Ineffective formulas (formulas without computable realisers) can be given an approximating realiser that works only when the memory has the correct information. The main result is that from a proof one can extract an agent that will be able to smartly build an approximation to the memory good enough to eventually produce a correct realiser. Finally, Alexander Miquel \cite{Miquel(2011)} has been working on extending Krivine's classical realizability with forcing, in the context of second-order arithmetic. This raises a few questions: 

\begin{itemize}
	\item[({\bf Q1})] What underlies the combination of realizability and forcing in general? Can forcing be combined with other functional interpretations, e.g. Diller-Nahm? What benefits could that bring?
	\item[({\bf Q2})] As with Goodman's interpretation, could the learning-based realizability be decomposed into a standard realizability interpretation followed by some variant of forcing? 
\end{itemize}

%%%%%%%%%%%%%%%%%%%%%%%%%%%%%%%%%%%%%%%%
\subsection{Bounded-like interpretations}
%%%%%%%%%%%%%%%%%%%%%%%%%%%%%%%%%%%%%%%%

Very recently \cite{BBS(2012)} variants of modified realizability and the \emph{dialectica} interpretation have been proposed which apply to proofs in \emph{nonstandard arithmetic}. The main feature of the interpretation is to extract from a proof of an existential statement a finite set of candidate witnesses (as in Herbrand's theorem), rather than a precise witness. The authors show that finite sets are the appropriate way to interpret existential \emph{standard} quantifiers, while unrestricted existential quantifiers are interpreted uniformly (as in \cite{Berger(05c)} and \cite{Krivine(2003)}). 

Also recently, so-called \emph{bounded} variants of the dialectica and modified realizability interpretations \cite{Ferreira(05A),Ferreira(05),Ferreira(06)} have been proposed which make use of the Howard/Bezem strong majorizability relation but in a more embedded way than Kohlenbach's monotone interpretation. The original motivation was to extend functional interpretations to deal with ineffective principles in analysis such as weak K\"onig's lemma even over weak fragments of analysis. The bounded modified realizability was then extended into a \emph{confined} variant \cite{FO(2008)} which looks both for upper and lower bounds. There are striking similarities between the functional interpretation of non-standard arithmetic and the bounded and confined interpretations, as pointed out in \cite{BBS(2012)}. That raises the question:

\begin{itemize}
	\item[({\bf Q3})] What is the common structure behind these \emph{bounded-like} interpretations? In joint work with Gilda Ferreira \cite{FO(2012)} we have extended the unifying framework to deal with the bounded and confined interpretations, but unfortunately, this does not look to be general enough to include the non-standard arithmetic interpretation  \cite{BBS(2012)}, as they make crucial use of a new form of functional application.
\end{itemize}

%%%%%%%%%%%%%%%%%%%%%%%%%%%%%%%%%%%%%%%%
\subsection{Type-free functional interpretations}
%%%%%%%%%%%%%%%%%%%%%%%%%%%%%%%%%%%%%%%%

We have so far only been discussing Kreisel's version of realizability known as \emph{modified realizability}. The original realizability interpretation, however, due to Kleene \cite{Kleene(45)}, makes use of numbers (codes of Turing machines) as realizers, rather than functionals of higher type. The crucial difference is that not all codes $n$ define a total function $\kleene{n} \colon \NN \to \NN$. As such, the realizability of an implication $A \to B$ was originally defined as
\[
\begin{array}{lcl}
	\kreal{n}{(A \to B)} & \pdefin & \forall k (\kreal{k}{A} \to \defined{\kleene{n}(k)} \wedge \kreal{\kleene{n}(k)}{B}),
\end{array}
\]
so $\kleene{n}$ only needs to be defined on $k$ if $k$ is indeed a realizer\footnote{To appreciate the difference between Kleene number realizability and Kreisel's modified realizability it is enough to point out that the former is sound for the Markov principle whereas the later isn't. In fact, Kreisel developed modified realizability \cite{Kreisel(59)} precisely to show that the Markov principle is independent of intuitionistic arithmetic.} for $A$. Let us refer to Kleene's original notion of realizability as \emph{number realizability}.  It is clear that a \emph{relational} variant of number realizability also exists. For instance, the clause for implication would be: 
\[
\begin{array}{lcl}
	\rkreal{n}{k}{(A \to B)} & \pdefin & \forall m \rkreal{k_0}{m}{A} \to \defined{\kleene{n}(k_0)} \wedge \rkreal{\kleene{n}(k_0)}{k_1}{B}
\end{array}
\]
where $k_0$ and $k_1$ denote the first and second projections inverses of the standard coding $\NN \times \NN \to \NN$.
That raises the following questions:
\begin{itemize}
%	\item[({\bf Q1})] Is there a ``relational" version of \emph{number} realizability, as there is one for Kreisel modified realizability (described in Section \ref{newReal})? The extra clause $\defined{\kleene{n}(k)}$ in the conclusion of the interpretation of $A \to B$ makes it difficult to apply to same reasoning and carry universal quantifiers around as part of the interpretation. The questions as posed above is obviously not completely precise, as we have not formally defined the notion of a ``relational version" of an interpretation. Although a positive answer would be easily verified, providing a negative answer would first involve actually making the question precise.
	%
	\item[({\bf Q4})] Is there a \emph{number realizability} interpretation of affine logic? By that we mean an interpretation which works on numbers rather than functionals of finite type, and makes use of the fact that realizers might be partial. For instance, that might involve modifying the clause for $A \lto B$ in (\ref{basic-inter}) as
	\[
	\begin{array}{lcl}
	\uInter{A \lto B}{n}{k} & \pdefin & \uInter{A}{k_0}{\kleene{n_1}(k)} \lto \uInter{B}{\kleene{n_0}(k_0)}{k_1}.
	\end{array}
	\]
	But the question is when should we require that $\kleene{n_0}(k_0)$ and $\kleene{n_1}(k)$ be defined so as to obtain not only a sound interpretation but also possibly interpret new principles that are not interpreted by Kreisel's modified realizability? It seems none of the obvious choices work. But that of course does not rule out more comprehensive changes which could lead to a sound interpretation.
	\item[({\bf Q5})] Related to ({\bf Q4}), can one in general show that every natural (e.g.~modular) functional interpretation of intuitionistic logic can be extended to an interpretation of intuitionistic affine logic? And, even if this is not the case, is it always possible to relate functional interpretations in a similar way to the one done in Section \ref{hybrid}, perhaps using different parameters than the interpretation of $\bang A$?
	
% The ``natural" adjective above aims to rule out artificial examples that make explicit use of the global syntax of formulas and proofs. One way to formalise ``natural" could be to insist that the both the formula translation and the proof translation components of the interpretation should be modular. 
	
	\item[({\bf Q6})] Is there a ``number variant" of the other aforementioned interpretations? Beeson \cite{Beeson(1978)} has looked at the question for the dialectica interpretation, which he calls a \emph{type-free} dialectica. Beeson points out that there cannot be one for the actual dialectica interpretation, as it requires decidability of quantifier-free formulas whereas statements of the form $\defined{\kleene{n}(k)}$ are not decidable in general. He then suggests a type-free variant of the Diller-Nahm interpretation as
	\[
	\begin{array}{lcl}
	\uInter{A \to B}{n}{k} & \pdefin & \defined{\kleene{n_1}(k)} \wedge (\forall i \in \kleene{n_1}(k) \, \uInter{A}{k_0}{i} \to \defined{\kleene{n_0}(k_0)} \wedge \uInter{B}{\kleene{n_0}(k_0)}{k_1}).
	\end{array}
	\]
	In other words, he requires the counter-example functions to be total\footnote{I confess to not have been able to completely verify the soundness of Beeson's interpretation. The problem seems to appear in the interpretation of the cut rule ($A \to B$ and $B \to C$ implies $A \to C$) as the ``positive" witnesses for $A \to B$ need not be total, but that is used in building the ``negative" witnesses for $A \to C$, which should be total (cf. \cite{Beeson(1978)} middle of page 221).}, whereas the witnessing functions might be partial. Could this be relaxed? Could this be translated to the setting of affine logic? Would this lead to extra principles that go beyond those interpreted by the typed Diller-Nahm interpretation?
\end{itemize}

%%%%%%%%%%%%%%%%%%%%%%%%%%%%%%%%%%%%%%%%
\subsection{Short games versus long games}
%%%%%%%%%%%%%%%%%%%%%%%%%%%%%%%%%%%%%%%%

The use of games between two players to model non-classical logics started with the work of Lorentzen \cite{Felscher(2002)} where formulas were put in correspondence with debates/dialogues so that those provable in intuitionistic logic corresponded to dialogues in which the first player had a winning ``strategy". This idea was refined in the works of Blass \cite{Blass(1992)}, Abramsky \cite{Abramsky(1994)} and several others, and led to complete semantics for fragments of linear logic. 

The connections between games and the functional interpretations such as G\"odel's dialectica have been there from the start \cite{Scott(1968)}. In the final section 8 of \cite{Blass(1992)}, Blass discusses at great length how one can view de Paiva's \cite{dePaiva(1989B)} categorical formulation of the Diller-Nahm interpretation of linear logic as arising from Blass' game semantics. Blass' suggestion is that the functional interpretation of linear logic arises by considering short two-move games combined according to his rules but including ``Skolemisation" steps whenever it may be necessary to bring a long game into a two-move game. 

\begin{itemize}
	\item[({\bf Q7})] I feel that a better understanding of the differences between long games with concrete moves and the short games with higher-order moves is still lacking. Although Blass shows how one can think of the dialectica category as arising from his game semantics, it is well known that \emph{dialectica-like} games are useful to interpret extra principles that go beyond the interpreted logic such as the Markov principle, independence of premise and the axiom of choice. Blass long games, however, capture precisely some fragments of the logic providing a sound and complete semantics. 
	\item[({\bf Q8})] Related to ({\bf Q7}), can functional interpretations be used to build fully abstract models? Another question that would provide guidance towards this is: How does the functional interpretation of the \emph{propositional} fragment of linear logic relate to other models of linear logic such as proof nets, monoidal closed categories, coherent spaces and phase semantics? 	
	\item[({\bf Q9})] In the context of long games people have been able to fine tune the interpreted logic by restricting the kind of strategies one or both of the players is allowed to play (e.g. innocent \cite{Hyland(2000)}, fair, history-free). Not much in this direction has been done in the setting of functional interpretations, whereby one could consider restrictions on the class of realisers in order to avoid interpreting certain principles. It seems hard, however, to think of any restrictions that would make the interpretation not sound with respect to the axiom of choice, for instance, as its realiser is the identity. But one could consider other restrictions such as linear functionals, functionals of certain complexity, etc.
	\item[({\bf Q10})] Using the nomenclature of game theory \cite{Fudenberg(1991)}, the long games considered by Blass and Abramsky are said to be in \emph{extensive form}. Such games can be thought of as trees where each node in the tree is assigned one of the players and terminal nodes determine which player has won. Games in extensive form can be brought into a so-called \emph{normal form}, a matrix specifying for each given pair of strategies for the two players which of the two wins the game if they follow these strategies. Games in normal form can also be thought of as two-move games. The two-move game arising from a functional interpretation is obviously not going to be the same as the normal form of the given strategic Blass/Abramsky game. Two questions arise: What is the relation between these two different two-move games that come for the same logical formula $A$? Moreover, could the functional interpretation way of constructing two-move games have any relevance to game theory?  
\end{itemize}

%%%%%%%%%%%%%%%%%%%%%%%%%%%%%%%%%%%%%%%%
\subsection{Treading between linear and intuitionistic logic}
%%%%%%%%%%%%%%%%%%%%%%%%%%%%%%%%%%%%%%%%

We have seen that we can better understand and generalise an interpretation of intuitionistic logic by moving to the more general (and finer) setting of affine logic. There are, however, some interesting logics in between linear (no contraction) and intuitionistic (full contraction) logic. For instance, consider the following ``intuitionistic" version of {\L}ukasiewicz logic ($\ILuL$) obtained by adding to affine intuitionistic logic the contraction schema 
\begin{equation} \label{weak-contraction}
A \lto S_B A \otimes K_B A
\end{equation}
where $S_B A \pdefin B \lto A$ and $K_B A \pdefin (A \lto B) \lto B$. Note that (\ref{weak-contraction}) clearly follows from $A \lto A \otimes A$ since $A$ implies over affine logic both $S_B A$ and $K_B A$. We can obtain ``classical" {\L}ukasiewicz logic ($\CLuL$) by adding the double negation elimination $(A\lNeg)\lNeg \lto A$. If we denote by $\CL =$ classical logic, $\IL = $ intuitionistic logic, $\ILL =$ intuitionistic affine logic, and $\CLL =$ classical affine logic, the relation between these six logics is shown in the diagram below, where an arrow from $X$ to $Y$ means that $Y$ is an extension of $X$.

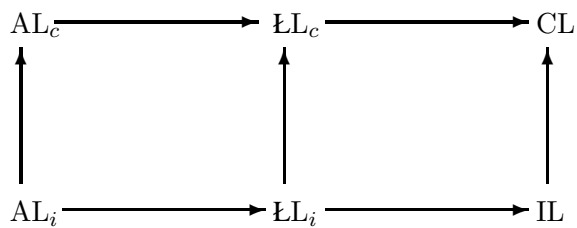
\begin{center}
	\setlength{\unitlength}{10mm}
	\begin{picture}(8,3)
		% border
		%\put(0.0,0.0){\line(1,0){8}}
		%\put(0.0,3.0){\line(1,0){8}}
		\thicklines
		% bottom row
		\put(0.5,0.0){$\ILL$}
		\put(1.2,0.15){\vector(1,0){2.7}}
		\put(4.0,0.0){$\ILuL$}
		\put(4.7,0.15){\vector(1,0){2.7}}
		\put(7.5,0.0){$\IL$}
		% top row
		\put(0.5,2.5){$\CLL$}
		\put(1.1,2.65){\vector(1,0){2.7}}
		\put(4.0,2.5){$\CLuL$}
		\put(4.7,2.65){\vector(1,0){2.7}}
		\put(7.5,2.5){$\CL$}
		% arrows going up
		\put(0.65,0.5){\vector(0,1){1.8}}
		\put(4.15,0.5){\vector(0,1){1.8}}
		\put(7.65,0.5){\vector(0,1){1.8}}
	\end{picture}
\end{center}

\begin{itemize}
	\item[({\bf Q11})] Since $\ILuL$ is a fragment of $\IL$, obviously any interpretation of $\IL$ also interprets $\ILuL$. The question, however, is whether one can make use of the fact that only limited contraction is available in $\ILuL$ and hence restrict the kind of functionals needed for the interpretation. For instance, which kind of minimal fragment of the simply-typed lambda calculus would be sufficient to provide a \emph{modified realizability} interpretation of $\ILuL$? This is related to {\bf (Q9)}. 

\end{itemize}

%%%%%%%%%%%%%%%%%%%%%%%%%%%%%%%%%%%%%%%%
\subsection{Endless possibilities?}
%%%%%%%%%%%%%%%%%%%%%%%%%%%%%%%%%%%%%%%%

The various functional interpretations discussed in Section \ref{hybrid} are only what one could call the ``classic" interpretations. As has been discussed in this Section \ref{future-work}, several other new and fascinating functional interpretations have been discovered recently. Beyond those already mentioned one also has:
\begin{itemize}
	\item Kohlenbach's \emph{monotone functional interpretations} \cite{Kohlenbach(A)}. These have been the cornerstone of the successful programme of \emph{proof mining} \cite{Kohlenbach(2008)}. It exploits a powerful combination of G\"odel's original dialectica interpretation with Howard's (or Bezem's) majorizability relation \cite{Bezem(85),Howard(73)}. 
	\item The \emph{Copenhagen interpretation} \cite{Biering(2008B)}. A variant of the dialectica interpretation where essentially in the interpretation of $A \to B$ the negative witnessing functional is allowed to ``give up" and not return a value. The original idea (apparently due to Martin Hyland) is that monads on types can quite often be lifted into an interpretation of (the comonad) $\bang A$. The Copenhagen interpretation carries this out for the monad $T X = X + 1$.
	\item \emph{Krivine classical realizability} \cite{Krivine(2003)}. Realizability interpretation of classical second order arithmetic, recently extended to countable choice. Krivine's realizability can be viewed as a combination of negative translation with a simpler intuitionistic realizability interpretation \cite{Oliva(2008A)}. 
\end{itemize}

We close with some final questions:

\begin{itemize}
	\item[({\bf Q12})] Is there a common structure behind all functional interpretations? What would be the appropriate way to define what functional interpretations are in general? 
	\item[({\bf Q13})] Functional interpretations of classical logic have all been shown to arise from an interpretation of intuitionistic logic combined with a negative translation. Can one show that this is always the case? 
\end{itemize}

\bibliographystyle{plain} 

\bibliography{../../../dblogic}

\end{document}